\documentclass[12pt]{article}
\usepackage{latexsym,amsmath,amssymb}
\usepackage{amsthm}
\usepackage{times}
\usepackage{amscd}
\usepackage{epsf}
\usepackage{graphicx}
\usepackage{graphics}

\newfont{\bb}{msbm10 at 12pt}
\newfont{\bbt}{msbm10 at 9pt}
\def\r{\hbox{\bb R}}

\def\b{\hbox{\bb B}}
\def\bt{\hbox{\bbt B}}
\def\rt{\hbox{\bbt R}}

\def\s{\hbox{\bb S}}
\def\st{\hbox{\bbt S}}
\def\disc{\mathcal{D}}
\def\mr{\mathbb{M}^2\times\mathbb{R}}

\def\m{\mathbb{M}^2}
\def\mt{\hbox{\bbt M}^2}

\newcommand{\parf}[2]{\dfrac{\partial #1 }{\partial #2}}

\newcommand{\norm}[1]{\left\Vert #1 \right\Vert}
\newcommand{\abs}[1]{\left\vert #1 \right\vert}
\newcommand{\set}[1]{\left\{#1\right\}}
\newcommand{\meta}[2]{\langle #1,#2 \rangle }
\newcommand{\eps}{\epsilon}

\newcommand{\To}{\longrightarrow }

\usepackage[latin1]{inputenc}
\usepackage {amsmath}
\usepackage {amsthm}
\usepackage{times}
\usepackage{amscd}
\usepackage{epsf}

\numberwithin{equation} {section}

\begin{document}

\theoremstyle{plain}\newtheorem{lemma}{Lemma}[section]
\theoremstyle{plain}\newtheorem{proposition}{Proposition}[section]
\theoremstyle{plain}\newtheorem{theorem}{Theorem}[section]
\theoremstyle{plain}\newtheorem{definition}{Definition}[section]
\theoremstyle{plain}\newtheorem{remark}{Remark}[section]
\theoremstyle{plain}\newtheorem{corollary}{Corollary}[section]

\begin{center}
\rule{15cm}{1.5pt} \vspace{.4cm}

{\Large \bf Fatou's Theorem and minimal graphs} \vspace{0.5cm}

{\large Jos$\acute{\text{e}}$ M. Espinar$\,^\dag$\footnote{The author is partially
supported by Spanish MEC-FEDER Grant MTM2007-65249, and Regional J. Andalucia Grants
P06-FQM-01642 and FQM325}, Harold Rosenberg$\,^\ddag$}\\
\vspace{0.4cm}\rule{15cm}{1.5pt}
\end{center}
\begin{flushright}
\today
\end{flushright}

\vspace{.4cm}

\noindent $\mbox{}^\dag$ Institut de Mathématiques, Universit$\acute{\text{e}}$
Paris VII, 175 Rue du Chevaleret, 75013 Paris, France; e-mail:
jespinar@ugr.es\vspace{0.2cm}

\noindent $\mbox{}^\ddag$ Instituto de Matematica Pura y Aplicada, 110 Estrada Dona
Castorina, Rio de Janeiro 22460-320, Brazil; e-mail: rosen@impa.br

\vspace{.3cm}

\begin{abstract}
In this paper we extend a recent result of Collin-Rosenberg ({\it a solution for the
minimal surface equation in the Euclidean disc has radial limits almost everywhere})
for a large class of differential operators in Divergence form. Also, we give an
alternative proof of Fatou's Theorem ({\it a harmonic function defined in the
Euclidean disc has radial limits almost everywhere}) even for harmonic functions
that are not bounded. Moreover, we construct an example (in the spirit of
\cite{CR2}) of a minimal graph in $\mr$, where $\m$ is a Hadamard surface, over a
geodesic disc which has finite radial limits in a mesure zero set.
\end{abstract}

\section{Introduction}

It is well known that a bounded harmonic function $u$ defined on the Euclidean disc
$D$ has radial limits almost everywhere (Fatou's Theorem \cite{F}). Moreover, the
radial limits can not be plus infinity for a positive measure set. For fixed $\theta
\in \s^{n-1}$, the radial limit $u(\theta)$ (if it exists) is defined as
\begin{equation*}
u(\theta) = \lim _{r\to 1}u(r,\theta),
\end{equation*}where we paramatrize the Euclidean disc in polar coordinates $(r,\theta)\in [0,1)\times \s
^1$.

In 1965, J. Nitsche \cite{N} asked if a Fatou Theorem is valid for the minimal
surface equation, i.e., {\it does a solution for the minimal surface equation in the
Euclidean disc have radial limits almost everywhere?} This question has been solved
recently by P. Collin and H. Rosenberg \cite{CR2}. Moreover, in the same paper
\cite{N}, J. Nitsche asked: {\it what is the largest set of $\theta$ for which a
minimal graph on $D$ may not have radial limits?} Again, this question was solved in
\cite{CR2} if one allows infinite radial limits. That is, they construct an example
of a minimal graph in the Euclidean disc with finite radial limits only on a set of
measure zero. In this example, the $+\infty$ radial limits (resp. $-\infty$) are
taken on a set of measure $\pi$ (resp. $\pi$).

The aim of this paper is to extend both results and give an alternative proof of
Fatou's Theorem for a more general situation. In Section \ref{Fatou}, we extend
Collin-Rosenberg's Theorem for a large class of differential operators in divergence
form (see Theorem \ref{Th:Fatou}). Also, we extend Fatou's Theorem even for
harmonic functions that are not bounded (see Theorem \ref{Th:Fatou2}). In
particular, as a consequence of this result, we obtain the classical Fatou Theorem
(see Corollary \ref{Cor:Fatou}). In Section \ref{Example}, we construct an example
of a minimal graph in $\mr $ over a geodesic disk $\disc \subset \m$ ($\m$ is a
Hadamard surface) for which the finite radial limits are of measure zero. Also, the
$+\infty$ radial limits (resp. $-\infty$) are taken on a set of measure $\pi$ (resp.
$\pi$).

\section{Fatou's Theorem}\label{Fatou}

Henceforth $(\b , g)$ denotes the $n-$dimensional unit open ball, i.e,
$$\b =\set{(r,\theta) \, ; \, 0\leq r < 1, \, \theta \in \s ^{n-1}},$$in polar coordinates
with respect to $g$, $g$ a $C^2-$Riemannian metric on $\b$. Define $G:= G(r,\theta)=
\sqrt{{\rm det}(g)}$. Moreover, we denote by $\nabla$ the Levi-Civita connection
associated to $g$ and by ${\rm div}_g$ its associated divergence operator. Also,
$L^1 (\b) $ denotes the set of integrable functions on $(\b , g)$.

Set $u \in C^2 (\b)-$function and $X_u$ be a $C^1(\b)-$vector field so that its
coordinates depend on $u$, its first derivatives and $C^1(\b)-$functions.

For fixed $\theta \in \s^{n-1}$, the radial limit ( if it exists) $u(\theta)$ is defined as
\begin{equation*}
u(\theta) = \lim _{r\to 1}u(r,\theta) .
\end{equation*}

\begin{theorem}\label{Th:Fatou}
Let $(\b , g, G , u, X_u)$ be as above. Assume that
\begin{itemize}
\item[a)] $\alpha \leq G(r , \theta) \leq \beta $ for all $(r,\theta)\in [0,1)\times \s
^{n-1}$, $\alpha$ and $\beta$ positive constants.
\item[b)] $|X_u| \leq M $ on $\b$, i.e., $X_u$ is bounded on $\b$.
\item[c)] $g(\nabla u , X_u) \geq \delta \, |\nabla u| + h$, where $\delta $ is a positive
constant and $h \in L^1 (\b)$.
\end{itemize}

Let $f \in L^1(\b)$. If $u$ is a solution of
\begin{equation*}
{\rm div}_g (X_u) \geq \, (\, \text{ or }\, \leq \, ) \,f \, \text{ on } \, \b,
\end{equation*}then $u$ has radial limits almost everywhere.
\end{theorem}
\begin{proof}
First, let us prove the case
$${\rm div}_g (X_u) \geq \,  \,f .$$

For $r<1$ fixed, set $\b (r)$ the $n-$dimensional open ball of radius $r$. Let $\eta
: \r \To (0 , +1)$ be a smooth function so that $0< \eta '(x)<1$ for all $x\in \r$.
Define $\psi := \eta \circ u $.

On the one hand, by direct computations and {\it item c)}, we have
\begin{equation*}
\begin{split}
{\rm div}_g (\psi \, X_u) &= \psi \, {\rm div}_g (X_u) + g(\nabla \psi , X_u) \geq
\psi \, f+ \eta ' \, g(\nabla u , X_u)\\[3mm]
  &\geq  \psi \, f + \eta ' \, \left( \delta \, |\nabla u| + h \right) = \delta \, \eta ' \, |\nabla u|
  + \left( \psi \, f + \eta ' \, h \right) \\
  &=\delta \,  |\nabla \psi |   + \left( \psi \, f + \eta ' \, h \right) ,
\end{split}
\end{equation*}thus
\begin{equation}\label{eqA}
\int _{\bt (r)}{\rm div}_g (\psi \, X_u) \geq \delta \int _{\bt (r)} |\nabla \psi| +
C
\end{equation}where $C$ is some constant. This follows since $f$ and $h$ are
$L^1-$functions on $\b$.

On the other hand, by Stokes' Theorem and {\it items a)} and ${\it b)}$, we obtain
for $r<1$ fixed
\begin{equation}\label{eqB}
\begin{split}
\int _{\bt (r)} {\rm div}_g (\psi \, X_u) &= \int_{\partial \bt (r)} \psi \, g(X_u
,\upsilon) \leq \int_{\partial \bt (r)} M \\
 &=  M \int _{\theta \in \st^{n-1}} G(r,\theta )d\theta \leq M \,\beta \int _{\theta \in \st ^{n-1}}\\
 &= M \, \beta \, \omega _{n-1} ,
\end{split}
\end{equation}where $\upsilon$ is the outer conormal to $\partial \b (r)$ and $\omega
_{n-1}$ is the volume of $\s ^{n-1}$.

So, from \eqref{eqA}, \eqref{eqB} and letting $r$ go to one, we conclude that
$|\nabla \psi|$ is integrable in $\b$, i.e.,
\begin{equation}\label{eqC}
\int _{\bt } |\nabla \psi | < + \infty
\end{equation}

Since $\parf{\psi}{r} \leq |\nabla \psi|$, we have from Fubini's Theorem and
\eqref{eqC}
\begin{equation*}
\int _{\bt} \parf{\psi}{r} = \int _{\theta \in \st
^{n-1}}\left(\int_0^1\parf{\psi}{r}G(r,\theta)\, dr \right) d\theta < \infty .
\end{equation*}

Thus, as $G(r,\theta)$ is bounded below by a positive constant, for almost all
$\theta \in \s ^{n-1}$,
\begin{equation*}
\lim _{r\to 1}\psi (r , \theta) - \psi (0,0) = \int _0 ^1 \parf{\psi}{r}(r, \theta)
dr < \infty,
\end{equation*}that is, $\psi $ has radial limits almost everywhere. Since $\psi = \eta \circ
u$, we conclude $u$ has radial limits almost everywhere (which may be $\pm \infty$).

For
$$ {\rm div}_g (X_u) \leq  f ,$$we just have to follow the above proof by changing
$\eta : \r \To (-1 ,0)$ so that $0< \eta ' (x) < 1$ for all $x \in \r$.
\end{proof}

As we pointed out in the Introduction, in the spirit of Theorem \ref{Th:Fatou}, we
can give an alternative proof of Fatou's Theorem even for harmonic function that are
not bounded, i.e.,

\begin{theorem}\label{Th:Fatou2}
Let $(\b , g, G , u)$ be as above. Assume that $\alpha \leq G(r , \theta) \leq \beta $
for all $(r,\theta)\in [0,1)\times \s ^{n-1}$, $\alpha$ and $\beta$ positive
constants. If $u$ is a solution of
\begin{equation*}
{\rm div}_g (\nabla u) = \, 0 \, \text{ on } \, \b,
\end{equation*}then $u$ has radial limits almost everywhere.

\end{theorem}
\begin{proof}
For $r<1$ fixed, set $\b(r)$ the $n-$dimensional open ball of radius $r$. Let $\eta
: \r \To (0,1)$ be a smooth function so that $0< \eta '(x)<1$ for all $x\in \r$.
Define
$$\phi := \eta \circ u , \mbox{ and } \psi := \dfrac{\phi}{\sqrt{1+\abs{\nabla u}^2}}. $$

On the one hand, by direct computations, we have
\begin{equation*}
\begin{split}
{\rm div}_g (\psi \, \nabla u) &= \psi \, {\rm div}_g (\nabla u) + \meta{\nabla \psi
}{\nabla u} = \meta{\nabla \psi }{ \nabla u}\\[3mm]
  &=  \eta ' \dfrac{\abs{\nabla u}^2}{\sqrt{1+\abs{\nabla u}^2}} -
  \phi \dfrac{\meta{\nabla \abs{\nabla u}^2}{\nabla u}}{2(1+\abs{\nabla u}^2)^{3/2}}\\[3mm]
  & \geq \eta ' \abs{\nabla u} - \dfrac{\eta '}{\sqrt{1+\abs{\nabla u}^2}}-
  \phi \dfrac{\meta{\nabla \abs{\nabla u}^2}{\nabla u}}{2(1+\abs{\nabla u}^2)^{3/2}} \\[3mm]
  &\geq \abs{\nabla \phi} -1- \dfrac{\meta{\nabla \abs{\nabla u}^2}{\nabla u}}{2(1+\abs{\nabla u}^2)^{3/2}},
\end{split}
\end{equation*}since
\begin{equation*}
\nabla \psi = \eta ' \dfrac{\nabla u}{\sqrt{1+\abs{\nabla u}^2}} - \phi
\dfrac{\nabla \abs{\nabla u}^2}{2(1+\abs{\nabla u}^2)^{3/2}}.
\end{equation*}

Let us first bound the term
$$ \abs{\int _{\bt (r)}\dfrac{\meta{\nabla \abs{\nabla u}^2}{\nabla u}}{2(1+\abs{\nabla u}^2)^{3/2}}} .$$

Set $Y_u := \dfrac{\nabla u}{\sqrt{1+\abs{\nabla u}^2}}$, then
\begin{equation*}
\begin{split}
{\rm div}_g(Y_u) &= \dfrac{1}{\sqrt{1+\abs{\nabla u}^2}}{\rm div}_g(\nabla u) -
\dfrac{\meta{\nabla \abs{\nabla u }^2}{\nabla u}}{2(1+\abs{\nabla u}^2)^{3/2}} \\
 &= - \dfrac{\meta{\nabla \abs{\nabla u }^2}{\nabla u}}{2(1+\abs{\nabla u}^2)^{3/2}}
\end{split}
\end{equation*}since ${\rm div}_g(\nabla u) = 0$. Applying Stoke's Theorem we
obtain
\begin{equation*}
\int _{\bt (r)} {\rm div}_g(Y_u) = \int _{\partial \bt (r)} \meta{Y_u}{\upsilon} ,
\end{equation*}that is
\begin{equation*}
\begin{split}
\abs{ \int _{\bt (r)} \dfrac{\meta{\nabla \abs{\nabla u }^2}{\nabla
u}}{2(1+\abs{\nabla u}^2)^{3/2}} } & = \abs{\int _{\partial \bt (r)}
\meta{Y_u}{\upsilon}} \leq \int _{\partial \bt (r)} \abs{Y_u} \\[3mm]
&= \int _{\partial \bt (r)} \dfrac{\abs{\nabla u}}{\sqrt{1+\abs{\nabla
 u}^2}} \leq \int _{\partial \bt (r)}  \leq C ,
\end{split}
\end{equation*}for some positive constant $C$.

Thus
\begin{equation}\label{eqA2}
\int _{\bt (r)}{\rm div} (\psi \, X_u) \geq \int _{\bt (r)} |\nabla \phi| - \tilde
C,
\end{equation}for some positive constant $\tilde C$.

On the other hand, by Stokes' Theorem we obtain for $r<1$ fixed
\begin{equation}\label{eqB2}
\begin{split}
\int _{\bt (r)} {\rm div} (\psi \, X_u) &= \int_{\partial \bt (r)} \psi \,
\meta{X_u}{\upsilon} \leq \int_{\partial \bt (r)} \frac{\abs{\nabla u}}{\sqrt{1+\abs{\nabla u}^2}}\\
 &\leq  \int _{\partial \bt (r)} \leq C',
\end{split}
\end{equation}where $\upsilon$ is the outer conormal to $\partial \b(r)$ and $C'$ is some positive constant.

So, from \eqref{eqA2}, \eqref{eqB2} and letting $r$ go to one, we conclude that
$|\nabla \phi|$ is integrable in $\b$, i.e.,
\begin{equation}\label{eqC2}
\int _{\bt } |\nabla \phi | < \infty
\end{equation}

Since $\parf{\phi}{r} \leq |\nabla \phi|$, we have from Fubini's Theorem and
\eqref{eqC2}
\begin{equation*}
\int _{\bt} \parf{\phi}{r} = \int _{\theta \in \st
^{1}}\left(\int_0^1\parf{\phi}{r}G(r,\theta)\, dr \right) d\theta < \infty .
\end{equation*}

Thus, as $G(r,\theta)$ is bounded below by a positive constant, for almost all
$\theta \in \s ^{n-1}$,
\begin{equation*}
\lim _{r\to 1}\phi (r , \theta) - \phi (0,0) = \int _0 ^1 \parf{\phi}{r}(r, \theta)
dr < \infty,
\end{equation*}that is, $\phi $ has radial limits almost everywhere. Since $\phi = \eta \circ
u$, we conclude $u$ has radial limits almost everywhere (which may be $\pm \infty$).
\end{proof}

Then, as a consequence

\begin{corollary}\label{Cor:Fatou}
Let $u$ be a harmonic function defined over the Euclidean disc. Then $u$ has radial
limits almost everywhere.
\end{corollary}

\subsection{Applications}

%

Moreover, we will see now how Theorem \ref{Th:Fatou} applies to get radial limits
almost everywhere for minimal graphs in ambient spaces besides $\r ^3$. We work here
in  Heisenberg space, but it is not hard to check that we could work with minimal
graphs in a more general submersion (see \cite{LR}).

First, we need to recall some definitions in  Heisenberg space (see \cite{ADR}).
The Heisenberg spaces are $\r ^3$ endowed with a one parameter family of metrics
indexed by bundle curvature by a real parameter $\tau \neq 0$. When we say the {\it
Heisenberg space}, we mean $\tau = 1/2$, and we denote it by $\mathcal{H}$.

In global exponential coordinates, $\mathcal{H}$ is  $\r ^3$ endowed with the metric
$$ g = (dx^2+ dy^2) + (\frac{1}{2}(y dx -xdy) +dz)^2.$$

The Heisenberg space is a Riemannian submersion $\pi : \mathcal{H} \To \r$ over the
standard flat Euclidean plane $\r ^2$ whose fibers are the vertical lines, i.e.,
they are the trajectories of a unit Killing vector field and hence geodesics.

Let $S_0 \subset \mathcal{H}$ be the surface whose points satisfy $z=0$. Let $D
\subset \r^2$ be the unit disc. Henceforth, we identify domains in $\r ^2$ with its
lift to $S _0$. The Killing graph of a function $u \in C^2 (D)$ is the surface
$$ \Sigma = \set{(x,y,u(x,y)) \, ; \, (x,y) \in D }. $$

Moreover, the minimal graph equation is
\begin{equation*}
{\rm div}_{\rt^2}(X_u) = 0,
\end{equation*}here ${\rm div}_{\rt^2}$ stands for the divergence operator in $\r ^2$
with the Euclidean metric $\meta{}{}$, and
$$ X_u := \frac{\alpha }{W} \partial _x + \frac{\beta}{W} \partial _y,$$where
\begin{equation*}
\alpha := \frac{y}{2}+u_x , \quad \beta : = -\frac{x}{2}+u_y,
\end{equation*}and
\begin{equation*}
W^2 = 1 + \alpha ^2 + \beta ^2 .
\end{equation*}

Thus, for verifying $u $ has radial limits almost everywhere (which may be $\pm
\infty$), we have to check conditions $a)$, $b)$ and $c)$. Item $a)$ is immediate
since we are working with the Euclidean metric.

Item $b)$ follows from
$$ \abs{X_u}^2 = \frac{\alpha ^2 + \beta ^2 }{1 + \alpha ^2 +\beta ^2} \leq 1 .$$

Now, we need to check Item $c)$. On one hand, using polar coordinates $x= r \cos
\theta$ and $y= r \sin \theta$, we have
\begin{equation*}
\begin{split}
W^2 &=1+ \alpha ^2 + \beta ^2 = 1 + u_x ^2 + u_y ^2 + (y u_x - xu_y ) +\frac{x^2 +
y^2}{4} \\
 &= 1 + |\nabla u|^2 + \meta{\nabla u}{(-y,x)}+ \frac{x^2 + y ^2}{4} \\
 & \geq 1+ |\nabla u|^2 - |\nabla u||(-y,x)|+ \frac{x^2+y^2}{4} \\
 &= 1+ |\nabla u |^2 - r |\nabla u| + \frac{r^2}{4}
\end{split}
\end{equation*}thus,
\begin{equation*}
W \geq \sqrt{ 1 + \left( |\nabla u| - \frac{r}{2}\right)^2 } \geq \abs{|\nabla u| -
\frac{r}{2}}.
\end{equation*}

We need a lower bound for $W$ in terms of $|\nabla u|$. To do so, we distinguish two
cases:

{\bf Case $|\nabla u| \leq 5/4$:} Since
$$ 1 - r \abs{\nabla u}+ \frac{r^2}{4} \geq 1- \frac{5r}{4} + \frac{r^2}{4} \geq 0 \mbox{ for all } r \leq 1,
$$we obtain
$$ W \geq \sqrt{\abs{\nabla u}^2 +1 - r \abs{\nabla u} + \frac{r^2}{4}} \geq \abs{\nabla u} .$$

{\bf Case $\abs{\nabla u} > 5/4$}: We already know that
$$ W \geq \abs{|\nabla u | - \frac{r}{2}} ,$$thus, for $\abs{\nabla u}> 5/4$, it is
easy to see that
$$ \abs{|\nabla u| - \frac{r}{2}} \geq \frac{3}{10} \abs{\nabla u} \mbox{ for all } r \leq 1 .$$

So, in any case, for $\delta = 3/10 >0$
\begin{equation}\label{boundW}
W \geq \delta \abs{\nabla u} .
\end{equation}

On the other hand,
\begin{equation*}
\begin{split}
\meta{\nabla u}{X_u} &= \frac{u_x ^2 + u_y ^2 + \frac{1}{2}(y u_x - x u_y)}{W} \\[4mm]
 &= \frac{1+ u_x ^2 + u_y ^2 + (y u_x - x u_y ) + \frac{x^2 +y^2}{4}}{W} -
 \frac{1 + \frac{1}{2}(y u_x - x u_y ) + \frac{x^2 +y^2}{4}}{W}\\[4mm]
 &= \frac{W^2}{W} + h = W + h \geq \delta \abs{\nabla u} + h ,
\end{split}
\end{equation*}where we have used \eqref{boundW} and $h$ denotes the $L^1(D)-$function
$$ h = - \frac{1 + \frac{1}{2}(y u_x - x u_y ) + \frac{x^2 +y^2}{4}}
{\sqrt{1 + u_x ^2 + u_y ^2 + (y u_x - xu_y ) +\frac{x^2 + y^2}{4}}} ,$$that is, Item
$c)$ is satisfied. So,

\begin{corollary}
A solution for the minimal surface equation in the Heisenberg space defined over a
disc has radial limits almost everywhere (which may be $\pm \infty$).
\end{corollary}

\section{An example in a Hadamard surface}\label{Example}

The aim of this Section is to construct an example of a minimal graph in $\mr $ over
a geodesic disk $\disc \subset \m$ ($\m$ is a Hadamard surface) for which the finite
radial limits are of measure zero.

We need to recall preliminary facts about graphs over a Hadamard surface (see
\cite{GR} for details). Henceforth, $\m$ denotes a simply connected with Gauss
curvature bounded above by a negative constant, i.e., $K_{\mt} \leq c < 0$.

Let $p_0 \in \m$ and $\disc $ be the the geodesic disk in $\m$ centered at $p _0$ of
radius one. Re-scaling in the metric, we can assume that
$$ {\rm max}\set{ K_{\mt} (p) \, ; \, \, p \in \overline{\disc}} = -1 . $$

From the Hessian Comparison Theorem (see e.g. \cite{JK}), $\partial \disc$ bounds a
strictly convex domain. We assume that $\partial \disc$ is smooth, otherwise we can
work in a smaller disc. We identify $\partial \disc = \s ^1$ and orient it
counter-clockwise.

We say that $\Gamma $ is an {\it admissible polygon} in $\disc$ if $\Gamma $ is a
Jordan curve in $\overline{\disc}$ which is a geodesic polygon with an even number
of sides and all the vertices in $\partial \disc$. We denote by $A_1,B_1,\ldots, A_k
,B_k$ the sides of $\Gamma $ which are oriented counter-clockwise. Recall that any two
sides can not intersect in $\disc$. Set $D$ the domain in $\disc$ bounded by
$\Gamma$. By $|A_i|$ (resp. $|B_j|$), we denote the length of such a geodesic arc.

\vspace{.5cm}

\begin{theorem}[\cite{P}]\label{ThJS}
Let $\Gamma \subset \m$ be a compact polygon with an even number of geodesic sides
$A_1,B_1,A_2,B_2,\ldots ,A_n,B_n$, in that order, and denote by $D$ the domain with
$\partial D = \Gamma$. The necessary and sufficient conditions for the existence of
a minimal graph $u$ on $D$, taking values $+\infty$ on each $A_i$, and $-\infty $ on
each $B_j$, are the two following conditions:
\begin{enumerate}
\item $\sum _{i=1}^n |A_i| = \sum _{i=1}^n |B_i|$,
\item for each inscribed polygon $P$ in $D$ (the vertices of $P$ are among the vertices of
$\Gamma$) $P \neq D$, one has the two inequalities:
$$ 2 a(P) < |P| \mbox{ and } 2 b(P) < |P|.$$
\end{enumerate}

Here $a(P)=\sum_{A_j \in P}|A_j|$, $b(P)=\sum _{B_j \in P}|B_j|$ and $|P|$ is the
perimeter of $P$.
\end{theorem}

%

\vspace{.5cm}

The construction of this example follows the steps in \cite[Section III]{CR2}, but
here we have to be more careful in the choice of the first {\it inscribed square}
and the {\it trapezoids}. We need to choose them as {\it symmetric} as possible.

\vspace{.5cm}

\begin{center}
\begin{figure}[!h]
\hspace{-1cm}\epsfysize=9cm \epsffile{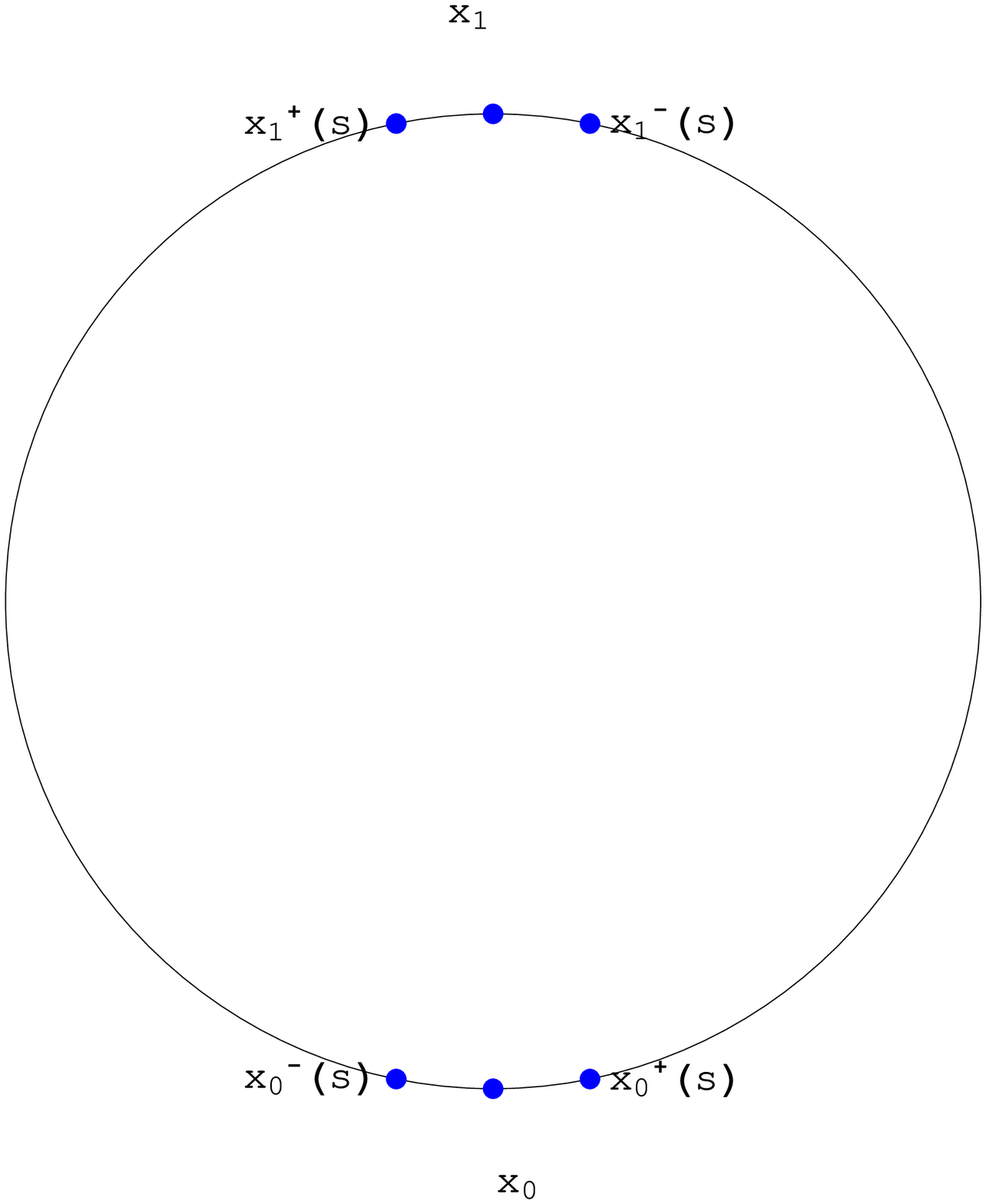}
\epsfysize=9cm \epsffile{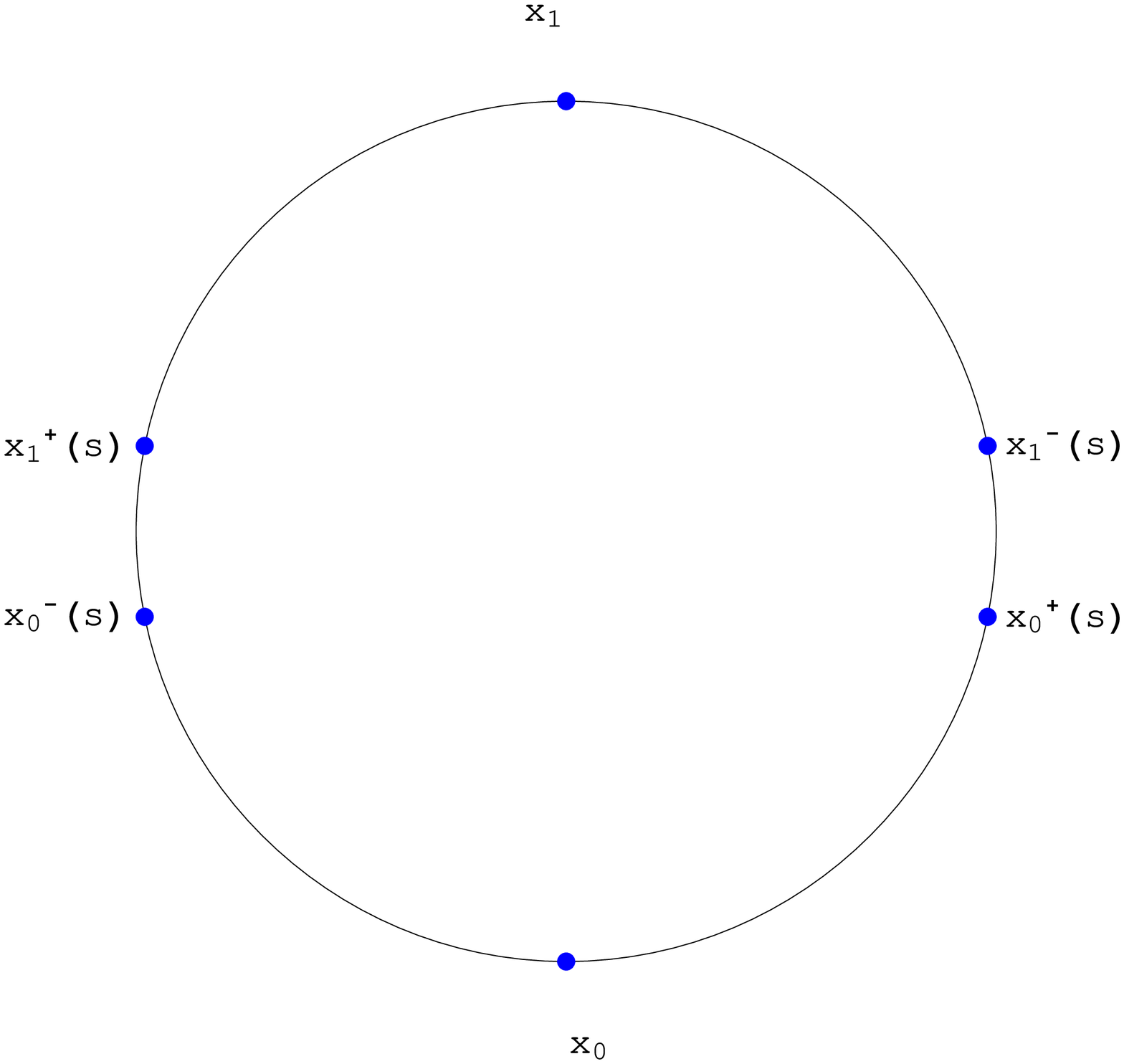}\\
\caption{We move the points along $\partial \disc$}\label{FigPto}
\end{figure}
\end{center}

Let us first explain how we take the {\it inscribed square}: Let $L = {\rm
length}(\partial \disc)$ and $\gamma (x_0 ,x_1)$ be the geodesic arc in $\disc$
joining $x_0 ,x_1 \in
\partial \disc $. Fix $x_0 \in \partial \disc $ and let $\alpha :\r / [0,L) \To \partial \disc
$ an arc-length parametrization of $\partial \disc $ (oriented count-clockwise). Set
$x_1 = \alpha(L/2)$. Consider $x_0 ^\pm (s) = \alpha (\pm s)$ and $x_1 ^± (s) =
\alpha (L/2 \pm s)$ for $0 \leq s \leq L/2$ (c.f. Figure \ref{FigPto}), and denote
\begin{eqnarray*}
B_1 (s) &=& \gamma ( x_0 ^+ (s) , x_1 ^- (s) )\\
A_1 (s) &=& \gamma ( x_1 ^- (s) , x_1 ^+ (s) )\\
B_2 (s) &=& \gamma ( x_1 ^+ (s) , x_0 ^- (s) )\\
A_2 (s) &=& \gamma ( x_0 ^- (s) , x_1 ^+ (s) ) .
\end{eqnarray*}

\begin{center}
\begin{figure}[!h]
\hspace{-1cm}\epsfysize=9cm \epsffile{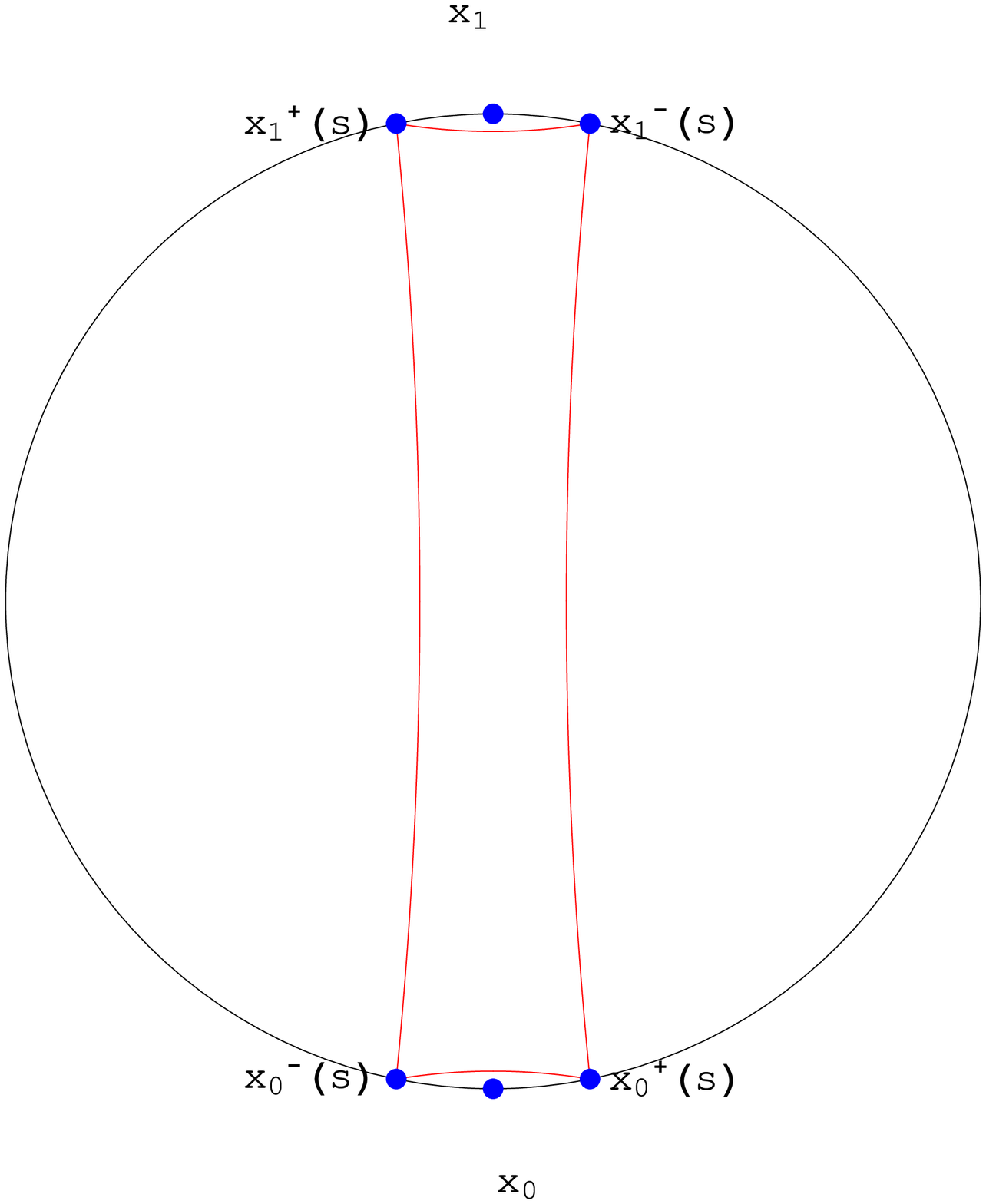}
\epsfysize=9cm \epsffile{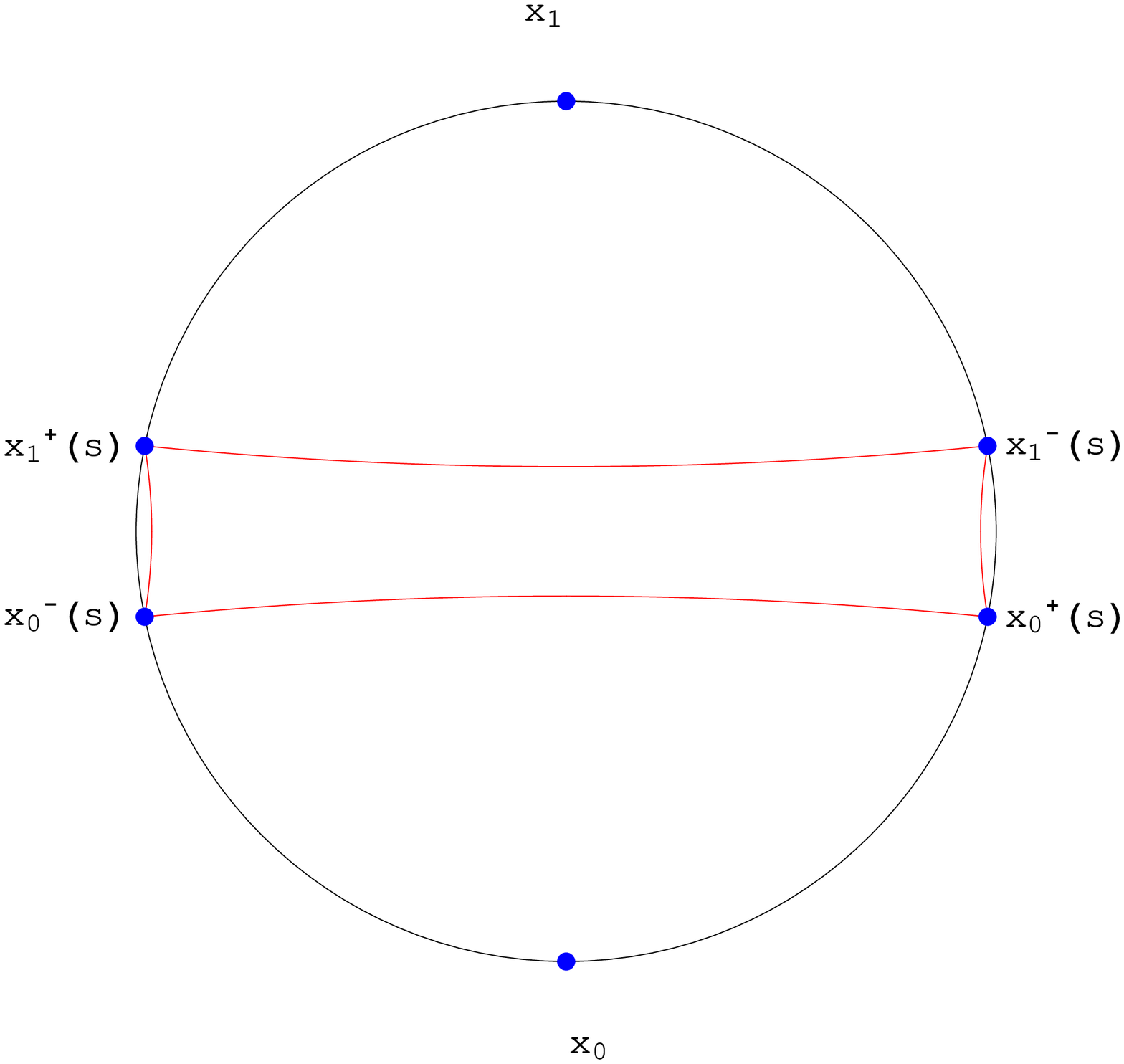}\\
\caption{How does the length change?}\label{FigPtoGeo}
\end{figure}
\end{center}

Hence (c.f. Figure \ref{FigPtoGeo}),
\begin{equation*}
\begin{matrix}
|A_1(s)|+ |A_2(s)| & > & |B_1(s)| + |B_2(s)| & \mbox{ for $s$ close to $0$.}\\
|A_1(s)|+ |A_2(s)| & < & |B_1(s)| + |B_2(s)| & \mbox{ for $s$ close to $L/2$.}
\end{matrix}
\end{equation*}

Thus, there exist $s_0 \in (0, L/2)$ so that
\begin{equation*}
|A_1(s_0)|+ |A_2(s_0)| = |B_1(s_0)| + |B_2(s_0)| .
\end{equation*}

So, given a fixed point $x_0 \in \partial \disc$, we have the existence of four
distinct points $p_1 = \alpha ( s_0)$, $p_2 = \alpha (L/2-s_0)$, $p_3 = \alpha (L/2
+s_0)$ and $p_4 = \alpha (-s_0)$ ordered counter-clockwise so that
\begin{equation*}
|A_1|+ |A_2| = |B_1| + |B_2| ,
\end{equation*}where
\begin{eqnarray*}
B_1 &=& \gamma ( p_1 , p_2 )\\
A_1 &=& \gamma ( p_2 , p_3 )\\
B_2 &=& \gamma ( p_3 , p_4 )\\
A_2 &=& \gamma ( p_4 , p_1 ) .
\end{eqnarray*}

\begin{figure}[!h]
\begin{center}
\epsfysize=9cm \epsffile{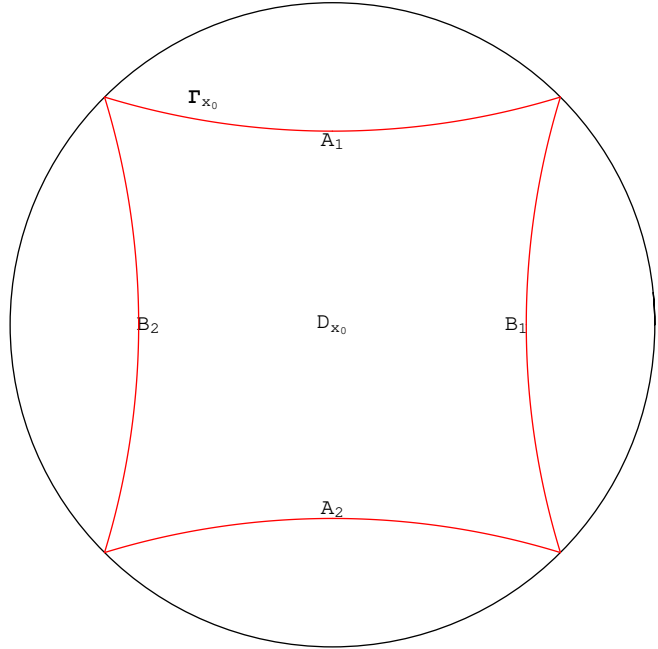}\\
\caption{Scherk domain}\label{FigScherkDom}
\end{center}
\end{figure}

In analogy with the Euclidean case \cite{CR2},

\begin{definition}\label{defsquare}
Fix a point $x_0 \in \partial \disc$, let $p_i$, $i=1,\ldots,4$ be the points
constructed above associated to $x_0 \in \disc$, then $\Gamma _{x_0}= A_1\cup B_1
\cup A_2 \cup A_3 $ is called the {\bf quadrilateral associated to} $x_0 \in \disc$
and it satisfies
\begin{equation*}
|A_1|+ |A_2| = |B_1| + |B_2| ,
\end{equation*}where
\begin{eqnarray*}
B_1 &=& \gamma ( p_1 , p_2 )\\
A_1 &=& \gamma ( p_2 , p_3 )\\
B_2 &=& \gamma ( p_3 , p_4 )\\
A_2 &=& \gamma ( p_4 , p_1 ) .
\end{eqnarray*}

Moreover, the interior domain $D_{x_0}$ bounded by $\Gamma _{x_0}$ is the {\bf
square inscribed associated to} $x_0 \in \disc$ (note that $D_{x_0}$ is a
topological disc), and $B_1$ is called the {\bf bottom side} (c.f. Figure
\ref{FigScherkDom}).
\end{definition}

Second, let us explain how to take the {\it regular trapezoids}: As above, fix $x_0
\in \partial \disc $ (from now on, $x_0$ will be fixed and we will omit it) and
parametrize $\partial \disc$ as $\alpha : \r / [0,L) \To
\partial \disc $. Let $0 \leq s_1 < s_2 < L$, or equivalently, two distinct and
ordered points $p_i =\alpha (s_i)\in \partial \disc$, $i=1,2$. The aim is to
construct a {\it trapezoid} in the region bounded by $\gamma (p_1, p_2)$ and
$\alpha([s_1 ,s_2])$. To do so, set $\bar{s} = \frac{s_1 + s_2}{2}$, i.e., $\bar{p}=
\alpha(\bar{s})$ is the mid-point. Define $\bar{p}^±(s) = \alpha (\bar{s} \pm s)$
for $0 \leq s \leq \bar{s}$.

Set
\begin{eqnarray*}
l_1 (s) &=& {\rm Length}\left(\gamma(p_1, \bar{p}^- (s))\right) \\
l_2 (s) &=& {\rm Length}\left(\gamma(\bar{p}^- (s), \bar{p}^+ (s))\right)\\
l_3 (s) &=& {\rm Length}\left(\gamma(\bar{p}^+ (s), p_2)\right)\\
l_4 (s) &=& {\rm Length}\left(\gamma(p_2, p_1)\right) .
\end{eqnarray*}

\begin{center}
\begin{figure}[!h]
\hspace{-1.9cm}\epsfysize=7cm \epsffile{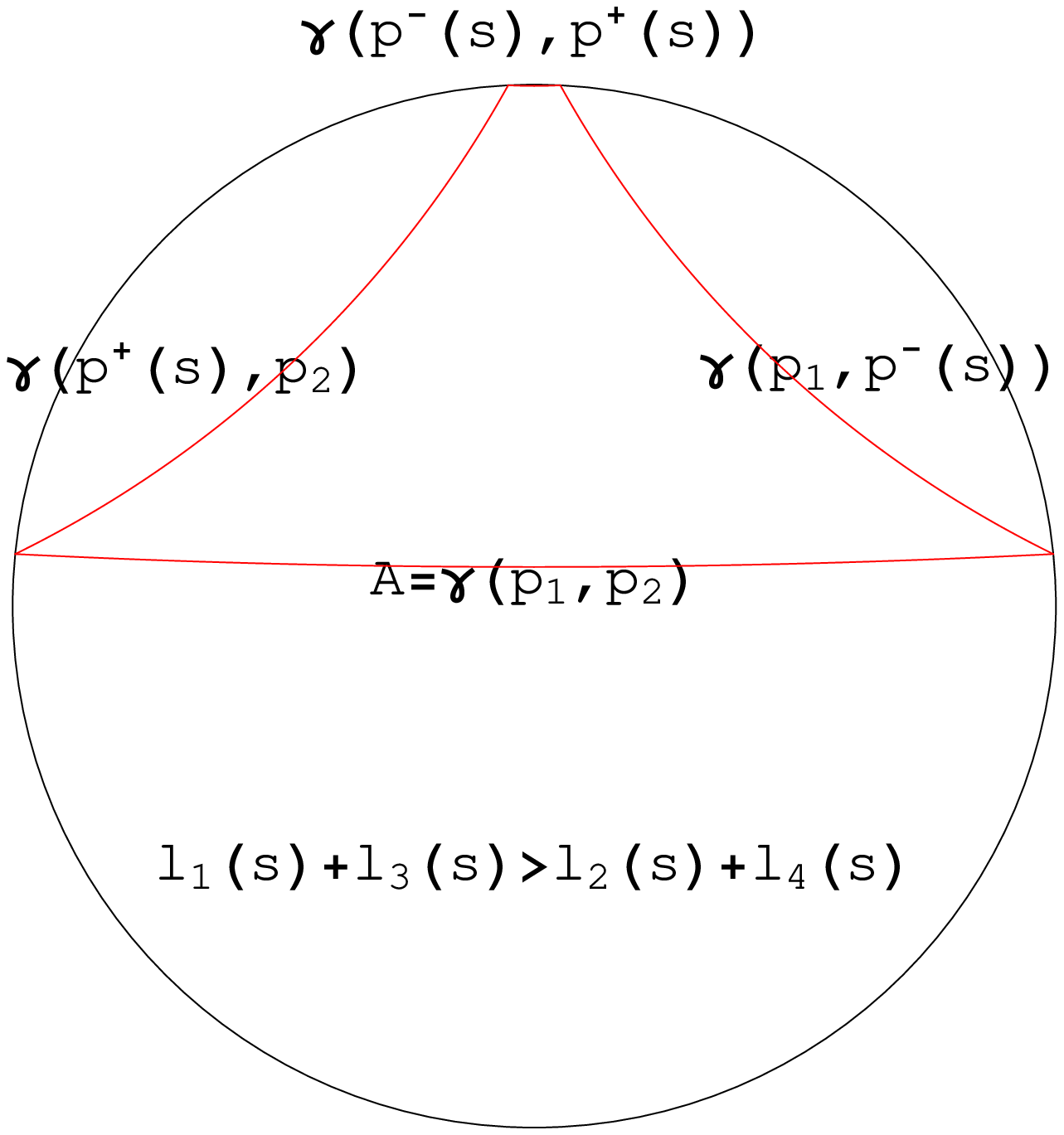}
\hspace{-1.1cm}\epsfysize=7cm \epsffile{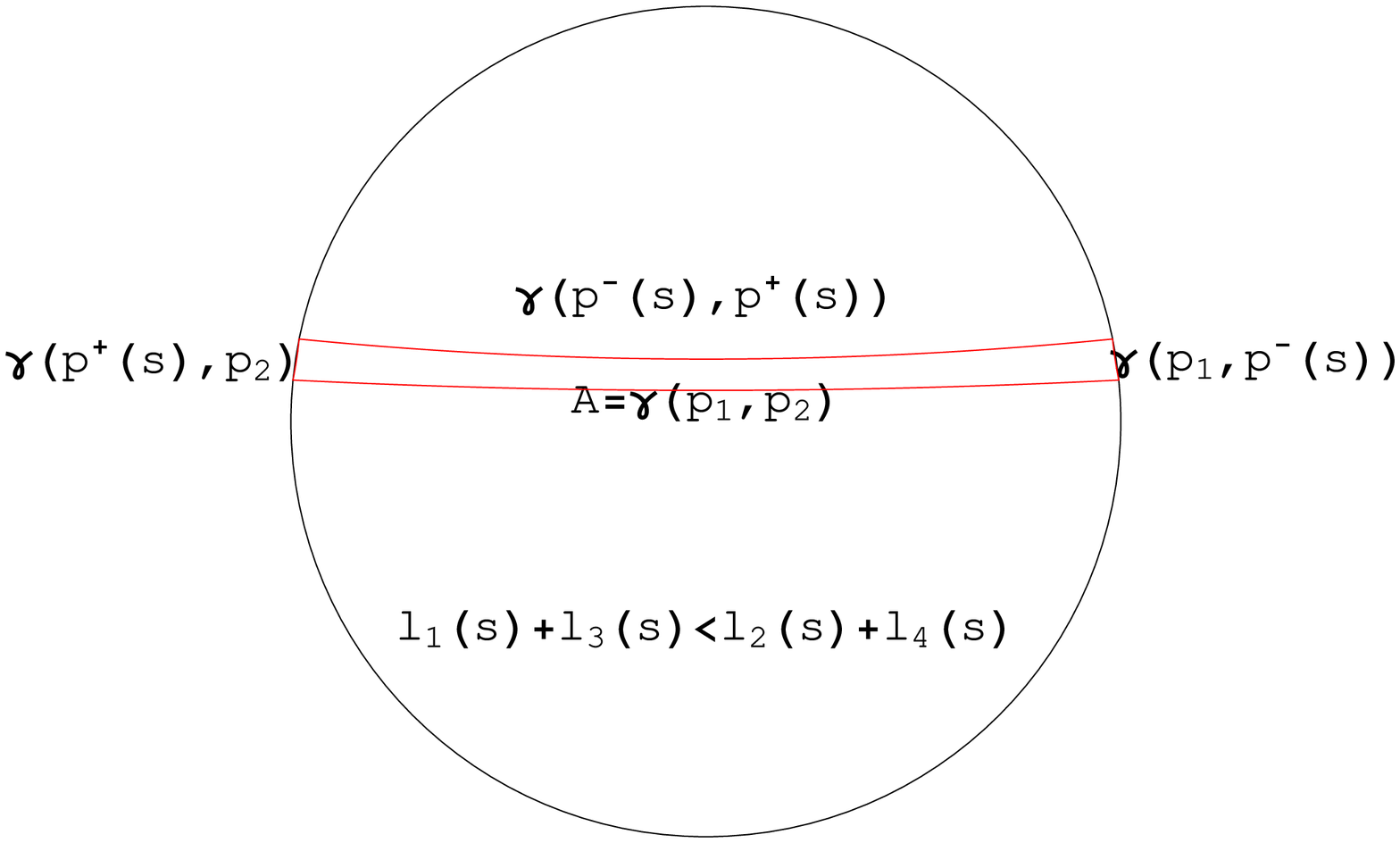}\\
\caption{How does the {\it trapezoid} vary?}\label{FigTrachange}
\end{figure}
\end{center}

Hence, for $s$ close to zero
\begin{equation*}
l_1(s) + l_3(s) > l_2(s) + l_4(s)
\end{equation*}by the Triangle Inequality, and for $s $ close to $\bar{s}$
\begin{equation*}
l_1(s) + l_3(s) < l_2(s) + l_4(s),
\end{equation*}since $l_1$ and $l_3$ go to zero and $l_4$ has positive length (c.f. Figure \ref{FigTrachange}).

Thus, there exists $s_0 \in (0, \bar{s})$ so that
\begin{equation*}
l_1(s_0) + l_3(s_0) = l_2(s_0) + l_4(s_0).
\end{equation*}

So, given a fixed point $x_0 \in \partial \disc$ and a geodesic arc $A := \gamma
(p_1, p_2)$ joining two (distinct and oriented) points in $\partial \disc$, we have
the existence of two distinct points $p^- = \alpha (\bar{s}-s_0)$ and $p^+ = \alpha
(\bar{s} +s_0)$ ordered count-clockwise so that
\begin{equation*}
l_1 + l_3 = l_2 + l_4 ,
\end{equation*}where
\begin{eqnarray*}
l_1 &=& {\rm Length}\left(\gamma(p_1, p^-)\right) \\
l_2 &=& {\rm Length}\left(\gamma(p^- , p^+ )\right)\\
l_3 &=& {\rm Length}\left(\gamma(p^+ , p_2)\right)\\
l_4 &=& {\rm Length}\left(\gamma(p_2, p_1)\right) .
\end{eqnarray*}

Moreover, the domain bounded by $\gamma(p_1, p^-) \cup \gamma(p^- , p^+ ) \cup
\gamma(p^+ , p_2) \cup \gamma (p_1 , p_2)$ is a topological disc.

\begin{figure}[!h]
\begin{center}
\epsfysize=5.5cm \epsffile{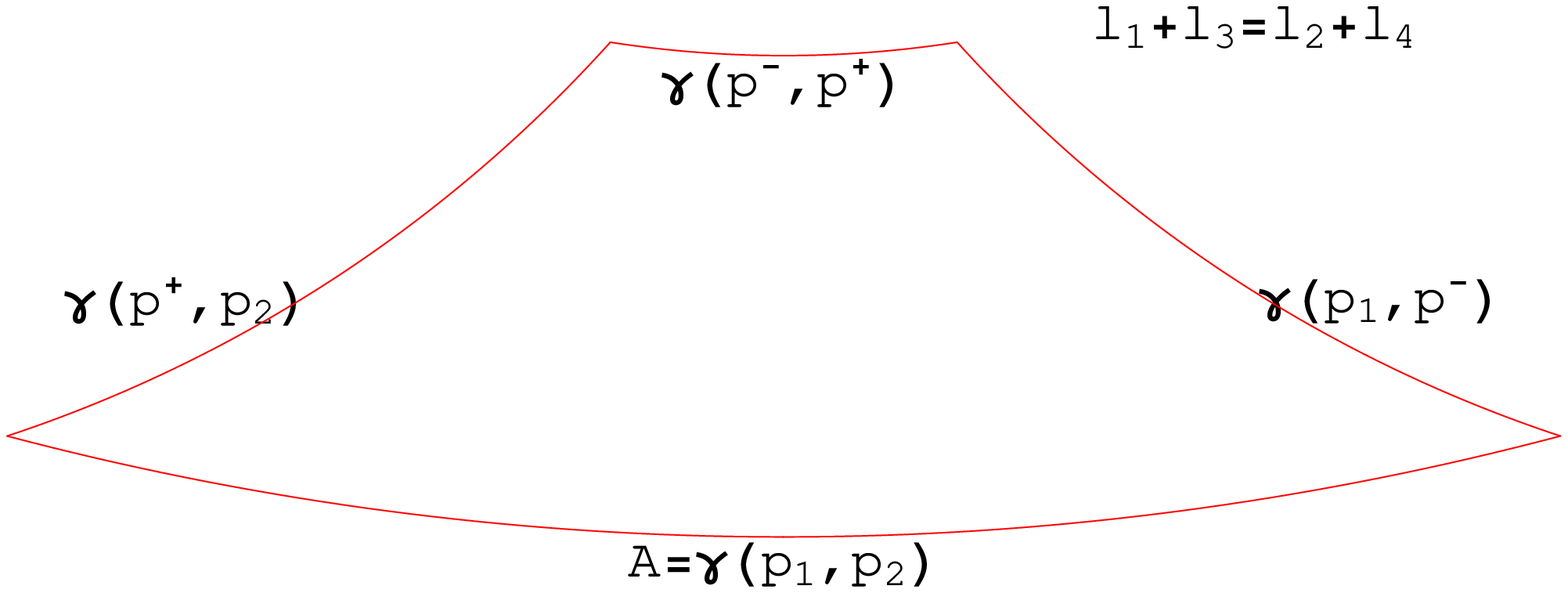}
\epsfysize=7.5cm \epsffile{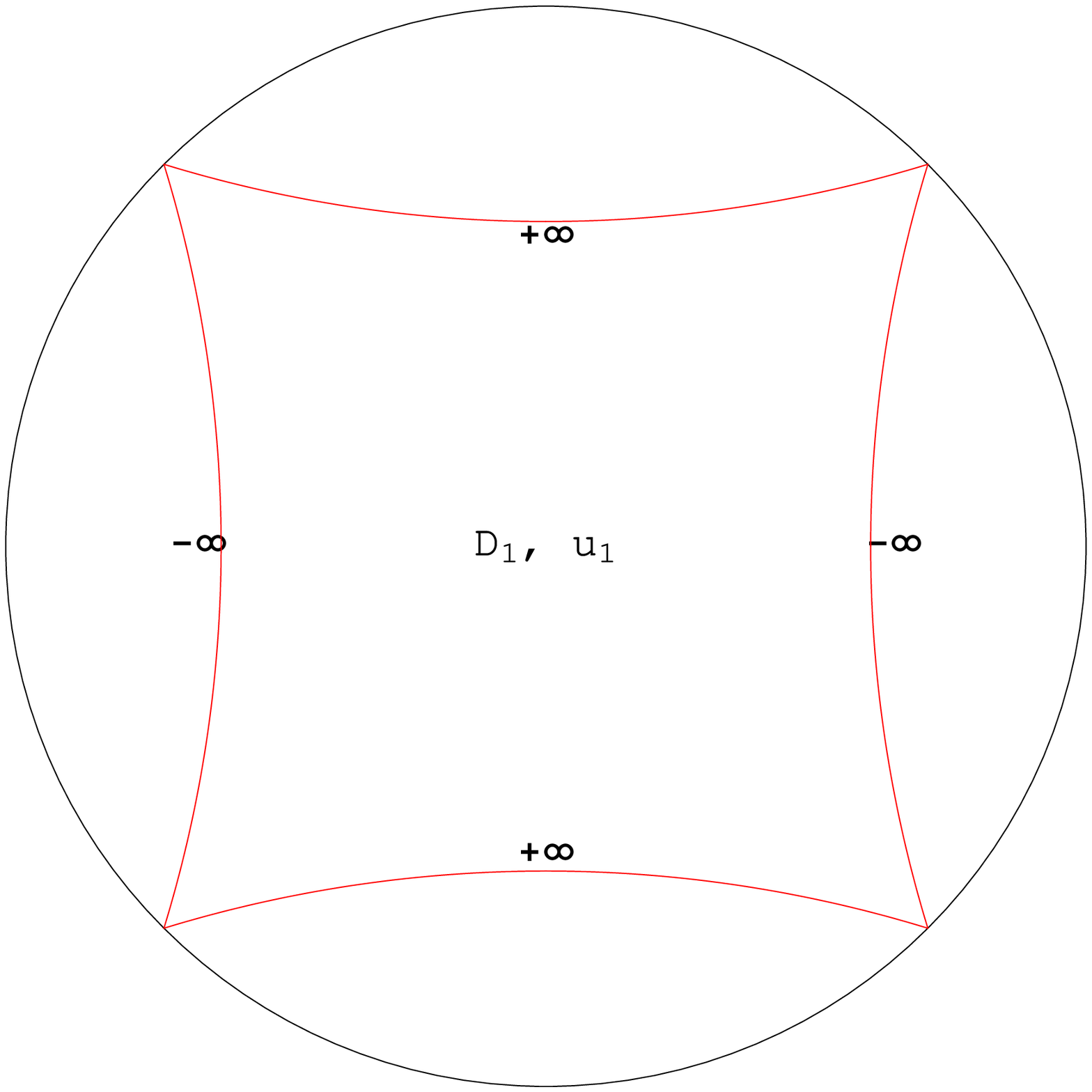}\\
\caption{(Left) Regular Trapezoid}\label{FigRegTra} \caption{(Right) First Scherk
domain}\label{FigD1}
\end{center}
\end{figure}


Again, in analogy with the Euclidean case,

\begin{definition}\label{deftrape}
$E = \gamma(p_1, p^-) \cup \gamma(p^- , p^+ ) \cup \gamma(p^+ , p_2) \cup \gamma
(p_1 , p_2)$ is called the {\bf regular trapezoid associated to the side} $A$, here
$A= \gamma(p_1 , p_2)$ (and, of course, once we have fixed a point $x_0 \in \partial
\disc$), and $p^±$ are given by the above construction (c.f. Figure
\ref{FigRegTra}).
\end{definition}

Now, we can begin the example. We only highlight the main steps in the construction
since, in essence, it is as in \cite[Section III]{CR2}.



Fix $x_0 \in \partial \disc$ and let $D_1$ the inscribed quadrilateral associated to
$x_0$ and $\Gamma _1 = \partial D_1$ (see Definition \ref{defsquare}). We label $A_1
, B_1 , A_2 , B_2$ the sides of $\Gamma _1$ ordered count-clockwise, with $B_1$ the
bottom side. By construction, $D_1$ is a Scherk domain. One can check this fact
using the Triangle Inequality. From Theorem \ref{ThJS}, there is a minimal graph
$u_1$ in $D_1$ which is $+ \infty$ on the $A_i's$ sides and equals $-\infty$ on the
$B_i 's$ sides (c.f. Figure \ref{FigD1}).

\begin{figure}[!h]
\begin{center}
\epsfysize=9cm \epsffile{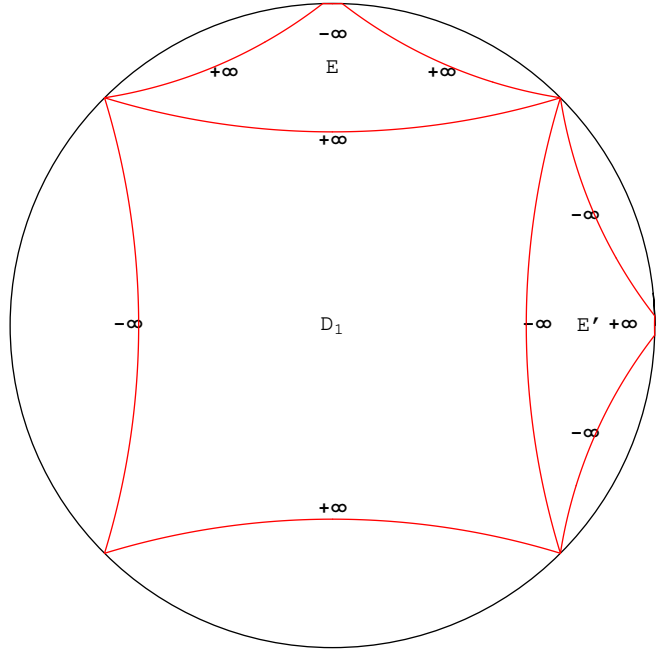}\\
\caption{Attaching trapezoids}\label{FigAttachTrap}
\end{center}
\end{figure}


Henceforth, we will attach regular trapezoids (see Definition \ref{deftrape}) to the
sides of the quadrilateral $\Gamma _1 $ in the following way. Let $E_1$ the regular
trapezoid associated to the side $A_1$, and $E' _1 $ the regular trapezoid
associated to the side $B_1$.

Consider the domain $D_2 = D_1 \cup E_1 \cup E' _1$, $\Gamma _2 = \partial D_2$.
This new domain does not  satisfy the second condition of Theorem \ref{ThJS}
,we only have to consider the inscribed polygon $E$ (c.f. Figure
\ref{FigAttachTrap}).

\begin{figure}[!h]
\begin{center}
\epsfysize=7.5cm \epsffile{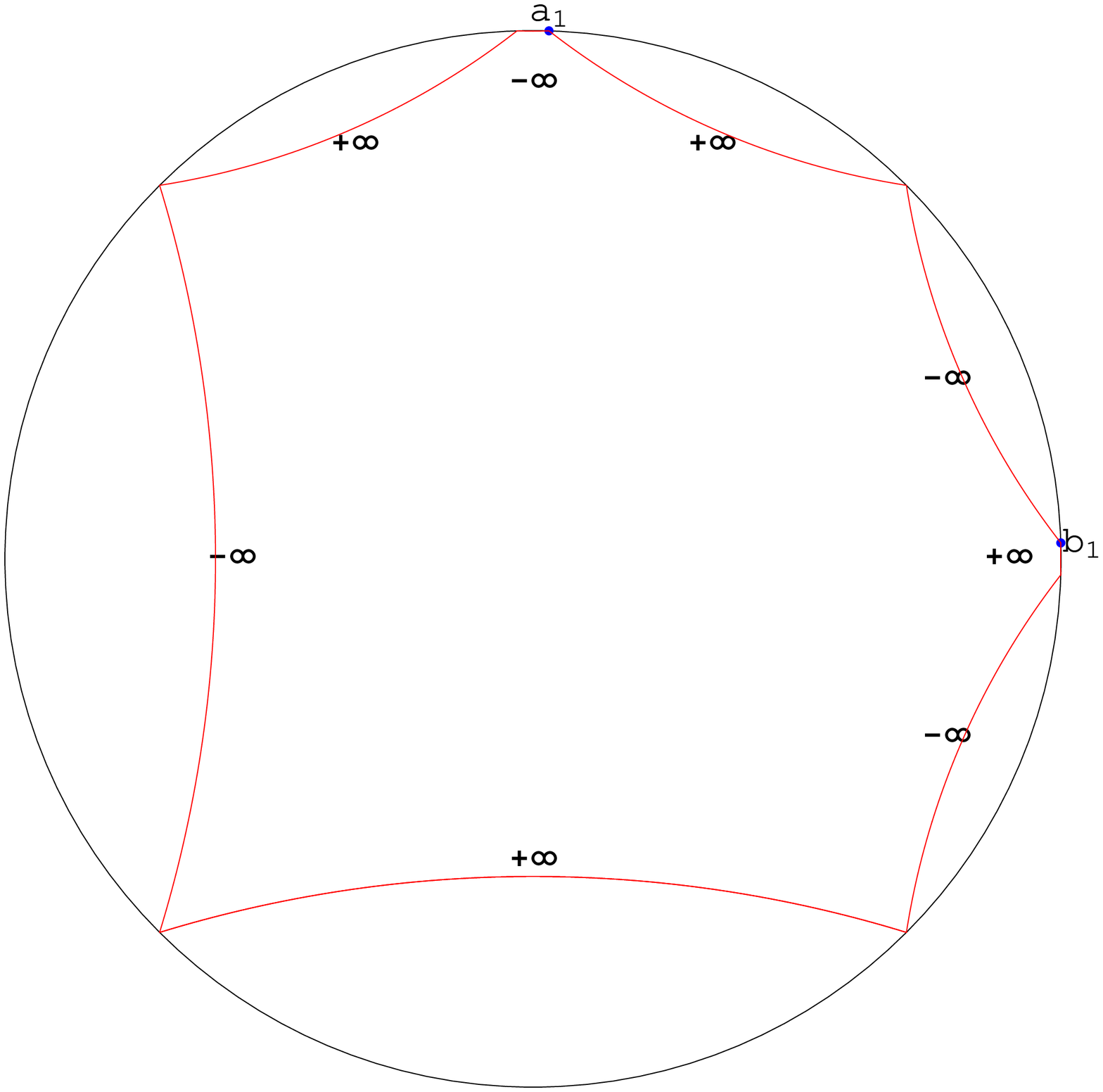}
\epsfysize=7.5cm \epsffile{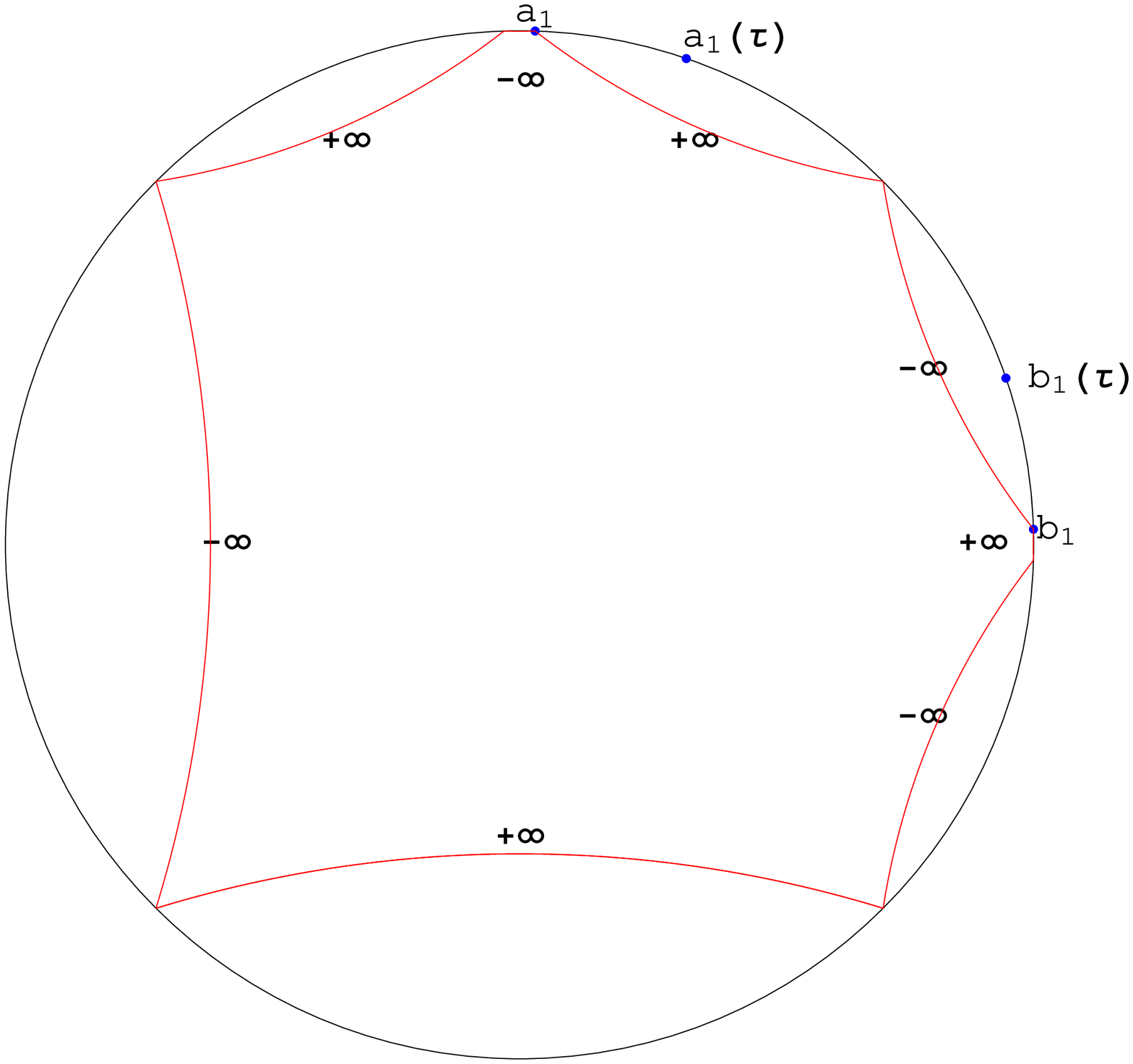}\\
\caption{Moving the vertex of the trapezoid}\label{FigPertPto}
\end{center}
\end{figure}


So, the next step is to perturb $D_2$ in such a way that it becomes an admissible
domain. Let $p $ be the common vertex of $E_1$ and $E ' _1$. Let $a_1 $ the closed
vertex of $E_1$ to $p$, and $b_1$ the closed vertex of $E'_1$ to $p$ (c.f. Figure
\ref{FigPertPto}).

One moves the vertex $a_1$ towards $b_1$ to a nearby point $a_1 (\tau)$ on $\partial
\disc$ (using the parametrization $\alpha : \r / [0,L) \To \partial \disc$ as we
have been done throughout this Section). And then one moves $b_1$ towards $a_1$ to a
nearby point $b_1 (\tau )$ on $\partial \disc $.

\begin{figure}[!h]
\begin{center}
\epsfysize=9cm \epsffile{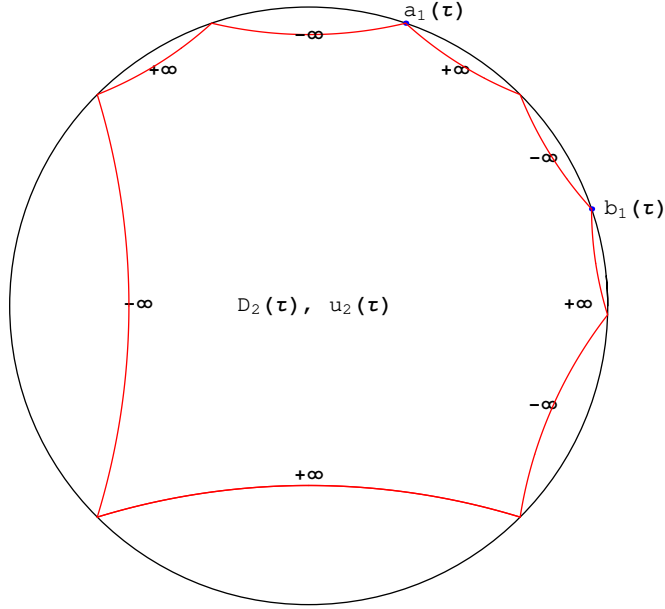}\\
\caption{Perturbed Scherk domain}\label{FigPertScherk}
\end{center}
\end{figure}

Let $\Gamma _2 (\tau )$ the inscribed polygon obtained by this perturbation, $E_1
(\tau)$ and $E'_1 (\tau)$ the perturbed regular trapezoids (c.f. Figure
\ref{FigPertScherk}). Thus, for $\tau > 0$ small, it is clear that:

\begin{itemize}
\item $\Gamma _2 (\tau)$  satisfies Condition 1 in Theorem \ref{ThJS}.
\item $2 a(E_1(\tau)) < |E_1(\tau)|$ and $2 b(E'_1(\tau)) < |E'_1 (\tau)|$.
\end{itemize}

Now, we state the following Lemma that establish how we extend the Scherk surface in
general.

\begin{lemma}\label{Lpert}
Let $u$ be a Scherk graph on a polygonal domain $D_1 = P(A_1, B_1, \ldots , A_k
,B_k)$, where the $A_i's$ and $B_i 's$ are the (geodesic) sides of $\partial D_1$ on
which $u$ takes values $+\infty$ and $-\infty$ respectively. Let $K$ be a compact
set in the interior of $D_1$. Let $D_2 = P(E_1, E'_1, A_2,B_2, \ldots, A_k,B_k)$ be
the polygonal domain $D_1$ to which we attach two regular trapezoids $E_1$ to the
side $A_1$ and $E'_1$ to the side $B_1$. Let $E_1 (\tau)$ and $E'_1 (\tau)$ be the
perturbed polygons as above. Then for all $\eps >0 $ there exists $\bar{\tau} >0 $
so that, for all $0< \tau \leq \bar{\tau}$, $v$ is a Scherk graph on $P(E_1 (\tau),
E'_1(\tau), A_2 , B_2 , \ldots, A_k , B_k)$ such that
\begin{equation}
\norm{u-v} _{C^2 (K)} \leq \eps .
\end{equation}
\end{lemma}
\begin{proof}
The proof of this Lemma relies on \cite[Section IV]{CR2} with the obvious
differences that we need to use the results for Scherk graphs over a domain in a
Hadamard surface stated in \cite{P} and \cite{GR}.
\end{proof}

Before we return to the construction, let us explain how we construct a {\it compact
domain associated to any Scherk domain}: Let $D = P(A_1, B_1 , \ldots , A_k,B_k)$ be
a Scherk domain in $\disc$ with vertex $\set{v_1, \ldots , v_{2k}} \in \partial
\disc$. Let $\beta _{v_i} :[0,1] \To \overline{\disc}$ denote the radial geodesic
starting at $p_0 \in \disc $ (the center of the disc $\disc$) and ending at $v_i \in
\partial \disc $. Note that any $\beta _{v_i}$ can not touch neither a $A_i$ side nor a
$B_i$ side expect at the vertex.

Set $r< 1$ and $p_i = \beta _{v_i} (r) \in \disc $ for $i = 1, \ldots , 2k$.
Consider the polygon
$$P = \bigcup _{i =1}^{2k-1}\gamma (p_{i},p_{i+i}) \cup \gamma
(p_{2k}, p_1) \subset D ,$$and let $K'$ be the closure of the domain bounded by $P$,
here $\gamma (p_i , p_{i+1})$ is the geodesic arc joining $p_i$ and $p_{i+1}$ in
$D$. Let $\mathcal{D}(p_i , 1-r)$ be geodesic disc centered at $p_i$ of radius $1-r$
for each $i = 1, \ldots ,2k$. Then,

\begin{definition}\label{defcomp}
For $r<1$ close to $1$, the {\bf compact domain associated to the Scherk domain} $D$
is given by
\begin{equation*}
K = K ' \setminus \bigcup _{i=1}^{2k} \mathcal{D}(p_i , 1-r) .
\end{equation*}
\end{definition}

\begin{figure}[!h]
\begin{center}
\epsfysize=7.5cm \epsffile{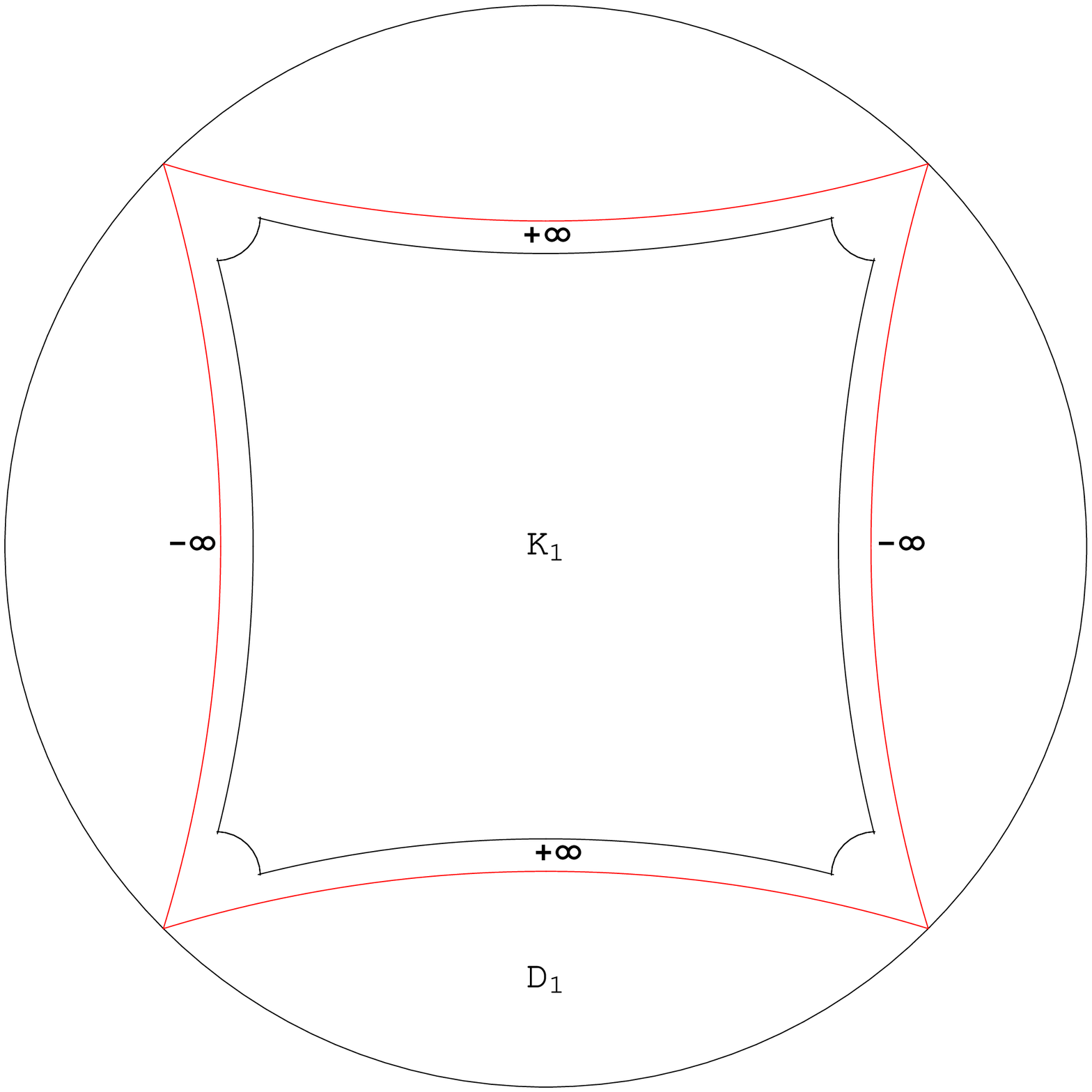}
\epsfysize=7.5cm \epsffile{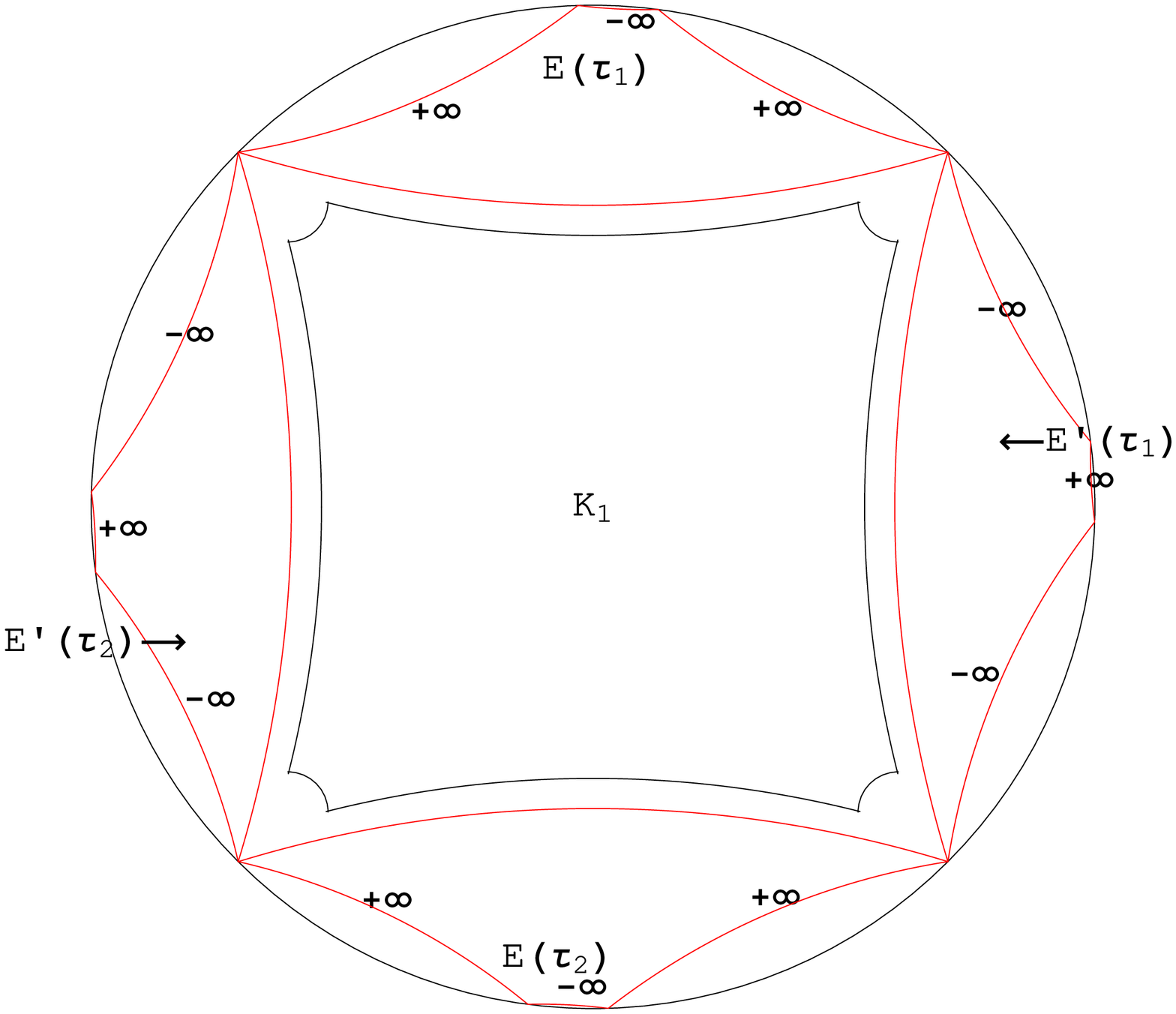} \\
\caption{(Left) Compact domain associated to the inscribed
quadrilateral}\label{FigK1} \caption{(Right) Attaching perturbed regular
trapezoids}\label{FigAllsides}
\end{center}
\end{figure}

Now, we continue with the construction. Let $D_1 =P(A_1,B_1,A_2,B_2)$ be the
inscribed square in $\disc$ (given in Definition \ref{defsquare}), and the Scherk
graph $u_1$ on $D_1$ which is $+\infty $ on the $A_i 's$ sides and $-\infty$ on the
$B_i 's$ sides. Let $K_1 $ be the compact domain associated to $D_1$ (see Definition
\ref{defcomp}). We choose $r_1 < 1$ close enough to one so that $u_1 >1$ on the
geodesic sides of $\partial K_1$ closer to the $A_i's$ sides and $u_1 < -1$ on the
geodesic sides of $\partial K_1$ closer to the $B_i 's$ sides (cf. Figure
\ref{FigK1}).

Next, we attach perturbed regular trapezoids to the sides $A_1$ and $B_1$, so from
Lemma \ref{Lpert}, for any $\eps _2 > 0$ there exists $\tau _2 >0 $ so that $D_2
(\tau ) = D_1 \cup E_1 (\tau) \cup E'_1(\tau)$ is a Scherk domain and $u_2 (\tau)$,
the Scherk graph defined on $D_2 (\tau)$, satisfy
\begin{equation*}
\norm{u_1 - u_2 (\tau)} _{C^2 (K_1)} \leq \eps _2 ,
\end{equation*}for all $0 < \tau \leq \tau _2$. Moreover, we can choose $u_2 (\tau
)$ so that $u_1 (p_0) = u_2 (\tau ) (p_0)$ (here $p_0$ is the center of $\disc$).
Then, choose $\eps _2 > 0 $ so that $u_2 (\tau ) > 1$ on the geodesic sides of
$\partial K_1$ closer to the $A_i's$ sides and $u_2 (\tau) < -1$ on the geodesic
sides of $\partial K_1$ closer to the $B_i 's$ sides.

\begin{figure}[!h]
\begin{center}
\epsfysize=7.5cm \epsffile{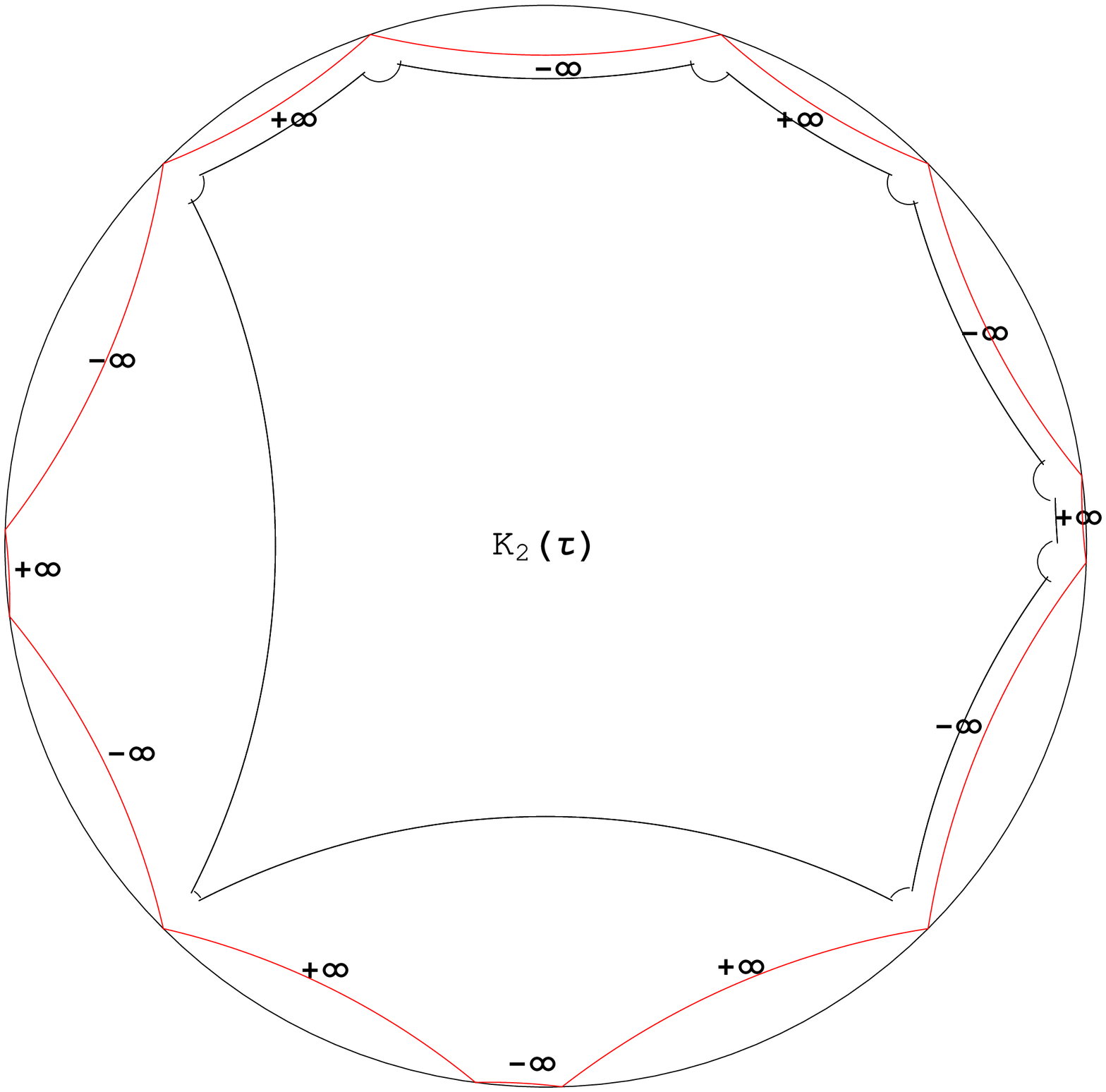} \epsfysize=7.5cm
\epsffile{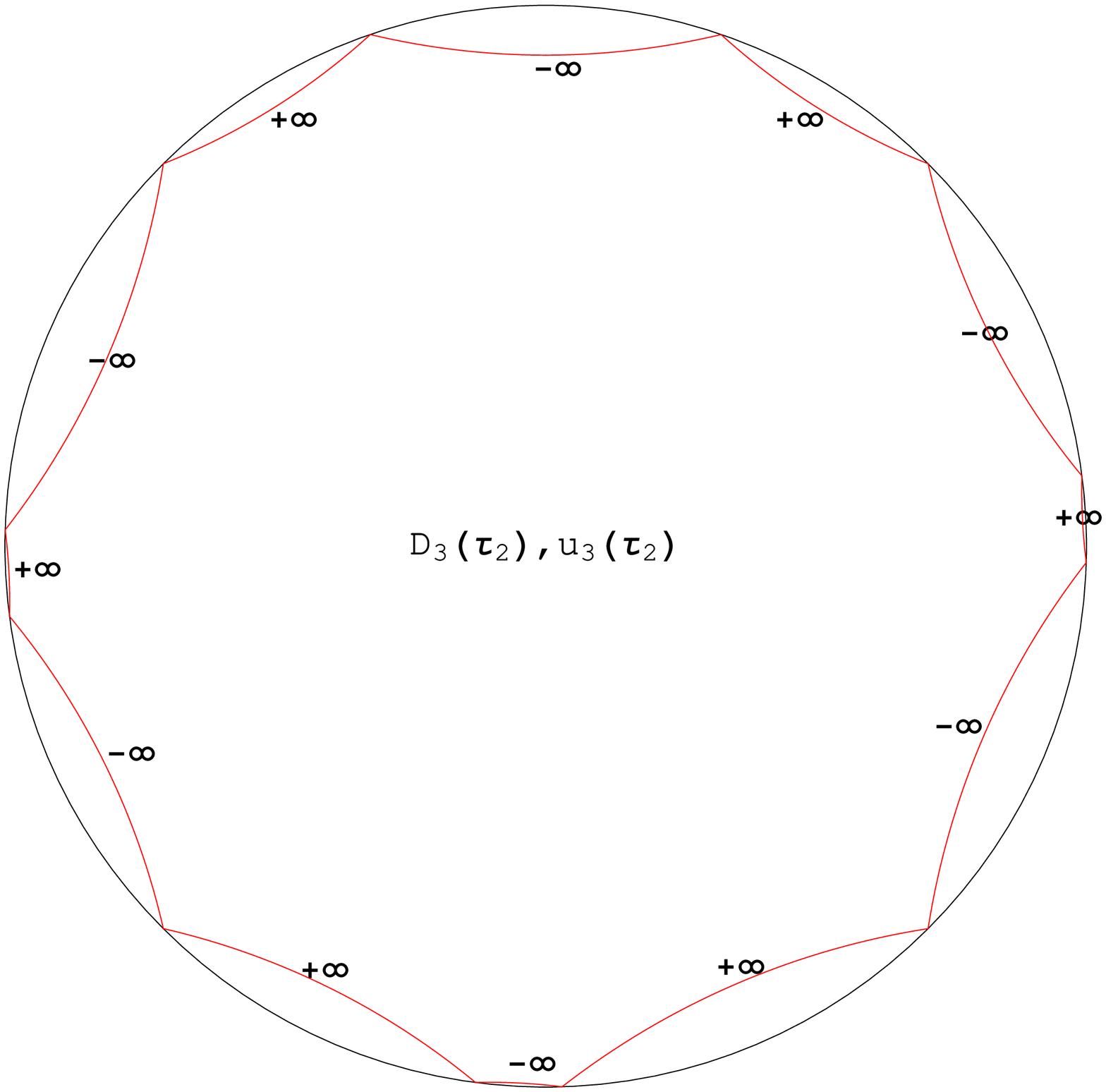}\\ \caption{(Left) Compact domain associated to $D_2
(\tau)$}\label{FigK2}  \caption{(Right) Choosing $u_3 (\tau)$}\label{FigD3}
\end{center}
\end{figure}

Let $K_2 (\tau)$ be the compact domain associated to the Scherk domain $D_2 (\tau)$.
Choose $r_2 <1 $ close enough to one (in the definition of $K_2 (\tau)$ given by
Definition \ref{defcomp}) so that, for $0 < \tau \leq \tau _2 $, $u_2 (\tau ) >2$ on
those geodesic sides of $\partial K_2 (\tau)$ parallel to the sides of $D_2 (\tau)$
where $u_2 (\tau ) = +\infty$, and $u_2 (\tau) <-2$ on the sides of $\partial K_2
(\tau)$ parallel to sides of $D_2 (\tau)$ where $u_2 (\tau ) = - \infty$ (cf. Figure
\ref{FigK2}).

Continue by constructing the Scherk domain $D_3 (\tau)$ by attaching perturbed
regular trapezoids (as above) to the sides $A_2$ and $B_2$ of $D_1$. We know, for
$\eps _3>0$, that there exist $\tau _3 >0$ so that if $0< \tau \leq \tau _3$ then
the Scherk graph $u_3(\tau)$ exists, $u_3(\tau)(p_0) = u_1(p_0)$ and
\begin{equation*}
\norm{u_3 (\tau) - u_2 (\tau)} _{C^2 (K_2(\tau))} \leq \eps _3 .
\end{equation*}

Moreover, choose $\eps _3 > 0 $ so that $u_3 (\tau ) > 3$ on the geodesic sides of
$\partial K_2 (\tau)$ closer to the $A_i's$ sides and $u_3 (\tau) < -3$ on the
geodesic sides of $\partial K_2 (\tau)$ closer to the $B_i 's$ sides (cf. Figure
\ref{FigD3}).


Now choose $\eps _ n \To 0$, $\tau _n \To 0$, $K_n (\tau _n)$ so that $K_n(\tau
_n)\subset  K_{n+1}(\tau_{n+1})$, $\bigcup_{n} K_n(\tau _n) = \disc$. Then the
$u_n(\tau _n)$ converge to a graph $u$ on $\disc$.

To see $u$ has the desired properties, we refer the reader to \cite[pages 13 and
14]{CR2} with the only difference that we need to use now Theorem \ref{Th:Fatou}.

\begin{remark}
The above construction can be carried out in a more general situation. Actually, if
we ask that
\begin{itemize}
\item The geodesic disc $\disc$ has strictly convex boundary.
\item There is a unique minimizing geodesic joining any
two points of the disc.
\end{itemize}

Then, we can extend the above example.
\end{remark}


\begin{thebibliography}{9}
\bibitem{ADR} L. Al\'{i}as, M. Dajczer, H. Rosenberg,
{\it The Dirichlet problem for constant mean curvature surfaces in Heisenberg
space}, Calc. Var. Partial Differential Equations, {\bf 30} (2007), nº 4, 513--522.
{\rm MR2332426}.

\bibitem{CR1} P. Collin, H. Rosenberg, {\it Construction of harmonic diffeomorphisms
and minimal graphs}. Preprint. {\rm Available at: arXiv:math.DG/0701547v1}.

\bibitem{CR2} P. Collin, H. Rosenberg, {\it Asymptotic values of minimal graphs ia disc}.
Preprint.

\bibitem{F} P. Fatou, {\it S\'{e}ries trigonom\'{e}triques et s\'{e}ries de Taylor},
Acta Math., {\bf 30} (1906), 335--400.

\bibitem{GR} J.A. Gálvez, H. Rosenberg, {\it Minimal surfaces and harmonic
diffeomorphisms from the complex plane onto a Hadamard surface}, preprint
(math.DG/0807.0997).

\bibitem{JK} L. Jorge, D. Koutroufiotis, {\it An estimate for the curvature of
bounded submanifolds}, Amer. J. Math., {\bf 103} (1981), no. 4, 711--725. {\rm
MR0623135}.

\bibitem{LR} C. Leandro, H. Rosenberg,
{\it Removable singularities for sections of Riemannian submersions of prescribed
mean curvature}. Preprint. {\rm Available at:
http://people.math.jussieu.fr/~rosen/Singularity.pdf}.

\bibitem{N} J. Nitsche, {\it On new results in the theory of minimal surfaces}, B. Amer. Math.
Soc., {\bf 71} (1965), 195--270. {\rm  MR0173993}.

\bibitem{P} A.L. Pinheiro, {\it A Jenkins-Serrin theorem in $\mr$}, To appear Bull.
Braz. Math. Soc.
\end{thebibliography}
\end{document}